\documentclass{amsart}
\usepackage{amssymb}
\usepackage{amsmath}
\usepackage{amsfonts}

\newtheorem{theorem}{Theorem}
\theoremstyle{plain}

\newtheorem{corollary}{Corollary}

\newtheorem{example}{Example}

\newtheorem{lemma}{Lemma}

\numberwithin{equation}{section}
\input{tcilatex}

\begin{document}
\title[M. Khan, Faisal and V.Amjid]{$\Gamma $-Abel-Grassmann's groupoids
characterized by their intuitionistic $\Gamma $-ideals}
\author{}
\maketitle

\begin{center}
$^{\ast }$\textbf{Madad Khan, Faisal and Venus Amjid}

\bigskip

\textbf{Department of Mathematics}

\textbf{COMSATS Institute of Information Technology}

\textbf{Abbottabad, Pakistan.}

\bigskip

$^{\ast }$\textbf{E-mail: madadmath@yahoo.com}

\bigskip

\textbf{In the memory of Professor Adam Khan}
\end{center}

\textbf{Abstract. }In this paper, we have discussed several properties of
intuitionistic fuzzy $\Gamma $-ideals of a $\Gamma $-AG-groupoid which is
the generalization of ideals in AG-groupoid We have characterized an
intra-regular $\Gamma $-AG$^{\ast \ast }$-groupoid in terms of
intuitionistic fuzzy $\Gamma $-left (right, two-sided) ideals,
intuitionistic fuzzy ($\Gamma $-generalized bi) $\Gamma $-bi-ideals,
intuitionistic fuzzy $\Gamma $-interior ideals and intuitionistic fuzzy $%
\Gamma $-quasi ideals. We have proved that the intuitionistic fuzzy $\Gamma $%
-left (right, interior, quasi) ideals coincide in an intra-regular $\Gamma $%
-AG$^{\ast \ast }$-groupoid. We have also shown that the set of
intuitionistic fuzzy $\Gamma $-two-sided ideals of an intra-regular $\Gamma $%
-AG$^{\ast \ast }$-groupoid forms a semilattice structure.

\textbf{Keywords. }$\Gamma $-AG-groupoids, intra-regular $\Gamma $-AG$^{\ast
\ast }$-groupoids and intuitionistic fuzzy $\Gamma $-ideals.

\begin{center}
\bigskip
\end{center}

\section{\textbf{Introduction}}

The concept of an Abel-Grassmann's groupoid (AG-groupoid) \cite{Kaz} was
first given by M. A. Kazim and M. Naseeruddin in $1972$ and they called it
left almost semigroup (LA-semigroup). P. Holgate call it simple invertive
groupoid \cite{hol}. An AG-groupoid $S$ is a groupoid having the left
invertive law $(ab)c=(cb)a$ for all $a,b,c\in S.$ An AG-groupoid is a
non-associative algebraic structure mid way between a groupoid and a
commutative semigroup. The left identity in an AG-groupoid if exists is
unique \cite{Mus3}. An AG-groupoid is non-associative and non-commutative
algebraic structure, nevertheless, it posses many interesting properties
which we usually find in associative and commutative algebraic structures.
An AG-groupoid with right identity becomes a commutative monoid \cite{Mus3}.
An AG-groupoid is basically the generalization of semigroup (see \cite{Kaz})
with wide range of applications in theory of flocks \cite{nas}.

Given a set $S$, a fuzzy subset of $S$ is an arbitrary mapping $%
f:S\rightarrow \lbrack 0,1]$ where $[0,1]$ is the unit segment of a real
line. This fundamental concept of fuzzy set was first given by Zadeh \cite%
{L.A.Zadeh} in $1965$. Fuzzy groups have been first considered by Rosenfeld 
\cite{A. Rosenfeld} and fuzzy semigroups by Kuroki \cite{N. Kuroki}.

Atanassov \cite{at}, introduced the concept of an intuitionistic fuzzy set.
Dengfeng and Chunfian \cite{8} introduced the concept of the degree of
similarity between intuitionistic fuzzy sets, which may be finite or
continuous, and gave corresponding proofs of these similarity measure and
discussed applications of the similarity measures between intuitionistic
fuzzy sets to pattern recognition problems. M. K. Sen introduced the concept
of $\Gamma $-semigroups in $1981.$ The non-associative $\Gamma $-AG-groupoid
is the generalization of an associative $\Gamma $-semigroup. In this paper,
we generalize several results in terms of intuitionistic fuzzy $\Gamma $%
-ideals.

Let $S$ and $\Gamma $ be any nonempty sets. If there exists a mapping $%
S\times \Gamma \times S\rightarrow S$ written as $(x,\alpha ,y)$ by $x\alpha
y,$ then $S$ is called a $\Gamma $-AG-groupoid if $x\alpha y\in S$ such that
the following $\Gamma $-left invertive law holds for\ all $x,y,z\in S$ and $%
\alpha ,\beta \in \Gamma $%
\begin{equation}
(x\alpha y)\beta z=(z\alpha y)\beta x.  \tag{$1$}
\end{equation}

A $\Gamma $-AG-groupoid also satisfies the $\Gamma $-medial law for\ all $%
w,x,y,z\in S$ and $\alpha ,\beta ,\gamma \in \Gamma $%
\begin{equation}
(w\alpha x)\beta (y\gamma z)=(w\alpha y)\beta (x\gamma z).  \tag{$2$}
\end{equation}

Note that if a $\Gamma $-AG-groupoid contains a left identity, then it
becomes an AG-groupoid with left identity.

A $\Gamma $-AG-groupoid is called a $\Gamma $-AG$^{\ast \ast }$-groupoid if
it satisfies the following law for\ all $x,y,z\in S$ and $\alpha ,\beta \in
\Gamma $%
\begin{equation}
x\alpha (y\beta z)=y\alpha (x\beta z).  \tag{$3$}
\end{equation}

A $\Gamma $-AG$^{\ast \ast }$-groupoid also satisfies the $\Gamma $%
-paramedial law for\ all $w,x,y,z\in S$ and $\alpha ,\beta ,\gamma \in
\Gamma $%
\begin{equation}
(w\alpha x)\beta (y\gamma z)=(z\alpha y)\beta (x\gamma w).  \tag{$4$}
\end{equation}

\section{Preliminaries}

Let $S$ be a $\Gamma $-AG-groupoid, a non-empty subset $A$ of $S$ is called
a $\Gamma $-AG-subgroupoid if $a\gamma b\in A$ for all $a$, $b\in A$ and $%
\gamma \in \Gamma $ or if$\ A\Gamma A\subseteq A.$

A subset $A$ of a $\Gamma $-AG-groupoid $S$ is called a $\Gamma $-left
(right) ideal of $S$ if $S\Gamma A\subseteq A$ $\left( A\Gamma S\subseteq
A\right) $ and $A$ is called a $\Gamma $-two-sided-ideal of $S$ if it is
both a $\Gamma $-left ideal and a $\Gamma $-right ideal.

A subset $A$ of a $\Gamma $-AG-groupoid $S$ is called a $\Gamma $%
-generalized bi-ideal of $S$ if $\left( A\Gamma S\right) \Gamma A\subseteq
A. $

A sub $\Gamma $-AG-groupoid $A$ of a $\Gamma $-AG-groupoid $S$ is called a $%
\Gamma $-bi-ideal of $S$ if $\left( A\Gamma S\right) \Gamma A\subseteq A$.

A subset $A$ of a $\Gamma $-AG-groupoid $S$ is called a $\Gamma $-interior
ideal of $S$ if $\left( S\Gamma A\right) \Gamma S\subseteq A.$

A subset $A$ of a $\Gamma $-AG-groupoid $S$ is called a $\Gamma $%
-quasi-ideal of $S$ if $S\Gamma A\cap A\Gamma S\subseteq A.$

\bigskip

A fuzzy subset $f$ is a class of objects with a grades of membership having
the form

\begin{center}
$f=\{(x,$ $f(x))/x\in S\}.$
\end{center}

An intuitionistic fuzzy set (briefly, $IFS$) $A$ in a non empty set $S$ is
an object having the form

\begin{center}
$A=\left\{ (x,\mu _{A}(x),\gamma _{A}(x))/x\in S\right\} .$
\end{center}

The functions $\mu _{A}:S\longrightarrow \lbrack 0,1]$ and $\gamma
_{A}:S\longrightarrow \lbrack 0,1]$ denote the degree of membership and the
degree of nonmembership respectively such that for all $x\in S,$ we have $%
0\leq \mu _{A}(x)+\gamma _{A}(x)\leq 1.$

For the sake of simplicity, we shall use the symbol $A=(\mu _{A},\gamma
_{A}) $ for an $IFS$ $A=\left\{ (x,\mu _{A}(x),\gamma _{A}(x))/x\in
S\right\} .$

Let $\delta =\left\{ (x,S_{\delta }(x),\Theta _{\delta }(x))/S_{\delta }(x)=1%
\text{ and }\Theta _{\delta }(x)=0/x\in S\right\} =(S_{\delta },\Theta
_{\delta })$ be an $IFS,$ then $\delta =(S_{\delta },\Theta _{\delta })$
will be carried out in operations with an $IFS$ $A=(\mu _{A},\gamma _{A})$
such that $S_{\delta }$ and $\Theta _{\delta }$ will be used in
collaboration with $\mu _{A}$ and $\gamma _{A}$ respectively.

An $IFS$ $A=(\mu _{A},\gamma _{A})$ of a $\Gamma $-AG-groupoid $S$ is called
an intuitionistic fuzzy $\Gamma $-AG-subgroupoid of $S$ if $\mu _{A}(x\alpha
y)\geq \mu _{A}(x)\wedge \mu _{A}(y)$ and $\gamma _{A}(x\alpha y)\leq \gamma
_{A}(x)\vee \gamma _{A}(y)$ for all $x$, $y\in S$ and $\alpha \in \Gamma .$

An $IFS$ $A=(\mu _{A},\gamma _{A})$ of a $\Gamma $-AG-groupoid $S$ is called
an intuitionistic fuzzy $\Gamma $-left ideal of $S$ if $\mu _{A}(x\alpha
y)\geq \mu _{A}(y)$ and $\gamma _{A}(x\alpha y)\leq \gamma _{A}(y)$ for all $%
x$, $y\in S$ and $\alpha \in \Gamma .$

An $IFS$ $A=(\mu _{A},\gamma _{A})$ of a $\Gamma $-AG-groupoid $S$ is called
an intuitionistic fuzzy $\Gamma $-right ideal of $S$ if $\mu _{A}(x\alpha
y)\geq \mu _{A}(x)$ and $\gamma _{A}(x\alpha y)\leq \gamma _{A}(x)$ for all $%
x,y\in S$ and $\alpha \in \Gamma .$

An $IFS$ $A=(\mu _{A},\gamma _{A})$ of a $\Gamma $-AG-groupoid $S$ is called
an intuitionistic fuzzy $\Gamma $-two-sided ideal of $S$ if it is both an
intuitionistic fuzzy $\Gamma $-left and an intuitionistic fuzzy $\Gamma $%
-right ideal of $S$.

An $IFS$ $A=(\mu _{A},\gamma _{A})$ of a $\Gamma $-AG-groupoid $S$ is called
an intuitionistic fuzzy $\Gamma $-generalized bi-ideal of $S$ if $\mu
_{A}((x\beta a)\gamma y)\geq \mu _{A}(x)\wedge \mu _{A}(y)$ and $\gamma
_{A}((x\beta a)\gamma y)\leq \gamma _{A}(x)\vee \gamma _{A}(y)$ for all $x$, 
$a$ and $y\in S$ and $\beta ,\gamma \in \Gamma $.

An intuitionistic fuzzy $\Gamma $-AG-subgroupoid $A=(\mu _{A},\gamma _{A})$
of a $\Gamma $-AG-groupoid $S$ is called an intuitionistic fuzzy $\Gamma $%
-bi-ideal of $S$ if $\mu _{A}((x\beta a)\gamma y)\geq \mu _{A}(x)\wedge \mu
_{A}(y)$ and $\gamma _{A}((x\beta a)\gamma y)\leq \gamma _{A}(x)\vee \gamma
_{A}(y)$ for all $x$, $a$ and $y\in S$ and $\beta ,\gamma \in \Gamma $.

An $IFS$ $A=(\mu _{A},\gamma _{A})$ of a $\Gamma $-AG-groupoid $S$ is called
an intuitionistic fuzzy $\Gamma $-interior ideal of $S$ if $\mu _{A}((x\beta
a)\gamma y)\geq \mu _{A}(a)$ and $\gamma _{A}((x\beta a)\gamma y)\leq \gamma
_{A}(a)$ for all $x$, $a$ and $y\in S$ and $\beta ,\gamma \in \Gamma $.

An $IFS$ $A=(\mu _{A},\gamma _{A})$ of a $\Gamma $-AG-groupoid $S$ is called
an intuitionistic fuzzy $\Gamma $-quasi ideal of $S$ if $(\mu _{A}\circ
_{\Gamma }S)\cap (S\circ _{\Gamma }\mu _{A})\subseteq \mu _{A}$ and $(\gamma
_{A}\circ _{\Gamma }S)\cup (S\circ _{\Gamma }\gamma _{A})\supseteq \gamma
_{A},$ that is, $(A\circ _{\Gamma }\delta )\cap (\delta \circ _{\Gamma
}A)\subseteq A.$

Let $S$ be a $\Gamma $-AG-groupoid and let $A_{I}=\{A/A\in S\},$ where $%
A=(\mu _{A},\gamma _{A})$ be any $IFS$ of $S,$ then $(A_{I},\circ _{\Gamma
}) $ satisfies $(1),$ $(2),$ $(3)$ and $(4)$.

An element $a$ of a $\Gamma $-AG-groupoid $S$ is called an intra-regular if
there exist $x,y\in S$ and $\alpha .\beta ,\gamma \in \Gamma $ such that $%
a=(x\alpha (a\beta a))\gamma y$ and $S$ is called an intra-regular if every
element of $S$ is an intra-regular.

\begin{example}
Let $S=\{1,2,3,4,5,6,7,8,9\}$. The following multiplication table shows that 
$S$ is an AG-groupoid and also an AG-band.
\end{example}

\begin{center}
\begin{tabular}{l|lllllllll}
. & $1$ & $2$ & $3$ & $4$ & $5$ & $6$ & $7$ & $8$ & $9$ \\ \hline
$1$ & $1$ & $4$ & $7$ & $3$ & $6$ & $8$ & $2$ & $9$ & $5$ \\ 
$2$ & $9$ & $2$ & $5$ & $7$ & $1$ & $4$ & $8$ & $6$ & $3$ \\ 
$3$ & $6$ & $8$ & $3$ & $5$ & $9$ & $2$ & $4$ & $1$ & $7$ \\ 
$4$ & $5$ & $9$ & $2$ & $4$ & $7$ & $1$ & $6$ & $3$ & $8$ \\ 
$5$ & $3$ & $6$ & $8$ & $2$ & $5$ & $9$ & $1$ & $7$ & $4$ \\ 
$6$ & $7$ & $1$ & $4$ & $8$ & $3$ & $6$ & $9$ & $5$ & $2$ \\ 
$7$ & $8$ & $3$ & $6$ & $9$ & $2$ & $5$ & $7$ & $4$ & $1$ \\ 
$8$ & $2$ & $5$ & $9$ & $1$ & $4$ & $7$ & $3$ & $8$ & $6$ \\ 
$9$ & $4$ & $7$ & $1$ & $6$ & $8$ & $3$ & $5$ & $2$ & $9$%
\end{tabular}
\end{center}

Clearly $S$ is non-commutative and non-associative because $2.3\neq 3.2$ and 
$(4.2).3\neq 4.(2.3).$

Let $\Gamma =\{\alpha ,\beta \}$ and define a mapping $S\times \Gamma \times
S\rightarrow S$ by $a\alpha b=(ab)^{2}$ and $a\beta b=a^{3}b^{2}$ for all $%
a,b\in S.$ Then it is easy to see that $S$ is a $\Gamma $-AG-groupoid and
also a $\Gamma $-AG-band. Note that $S$ is non-commutative and
non-associative because $9\alpha 1\neq 1\alpha 9$ and $(6\alpha 7)\beta
8\neq 6\alpha (7\beta 8).$

\section{Characterizations of an intra-regular $\Gamma $-AG$^{\ast \ast }$%
-groupoids in terms of intuitionistic fuzzy $\Gamma $-ideals}

\begin{example}
\label{ex}Let $S=\{1,2,3,4,5\}$ be an AG-groupoid with the following
multiplication table.
\end{example}

\begin{center}
\begin{tabular}{l|lllll}
. & $1$ & $2$ & $3$ & $4$ & $5$ \\ \hline
$1$ & $1$ & $1$ & $1$ & $1$ & $1$ \\ 
$2$ & $1$ & $2$ & $2$ & $2$ & $2$ \\ 
$3$ & $1$ & $2$ & $4$ & $5$ & $3$ \\ 
$4$ & $1$ & $2$ & $3$ & $4$ & $5$ \\ 
$5$ & $1$ & $2$ & $5$ & $3$ & $4$%
\end{tabular}
\end{center}

Let $\Gamma =\{\alpha \}$ and define a mapping $S\times \Gamma \times
S\rightarrow S$ by $x\alpha y=xy$ for all $x,y\in S,$ then clearly $S$ is a $%
\Gamma $-AG-groupoid.

It is easy to see that $S$ is an intra-regular. Define an $IFS$ $A=(\mu
_{A},\gamma _{A})$ of $S$ as follows: $\mu _{A}(1)=1$, $\mu _{A}(2)=$ $\mu
_{A}(3)=$ $\mu _{A}(4)=$ $\mu _{A}(5)=0,$ $\gamma _{A}(1)=0.3,$ $\gamma
_{A}(2)=0.4$ and $\gamma _{A}(3)=\gamma _{A}(4)=\gamma _{A}(5)=0.2,$ then
clearly $A=(\mu _{A},\gamma _{A})$ is an intuitionistic fuzzy $\Gamma $%
-two-sided ideal and also an intuitionistic fuzzy $\Gamma $-AG-subgroupoid
of $S$.

It is easy to observe that in an intra-regular $\Gamma $-AG-groupoid $S,$
the following holds%
\begin{equation}
S=S\Gamma S\text{.}  \tag{$5$}
\end{equation}

For an $IFS$ $A=(\mu _{A},\gamma _{A})$ of an AG-groupoid $S$ and $\alpha
\in (0,1],$ the set

$A_{\alpha }=\{x\in S:\mu _{A}(x)\geq \alpha ,$ $\gamma _{A}(x)\leq \alpha
\} $

is called an intuitionistic $\Gamma $-level cut of $A.$

\begin{theorem}
For a $\Gamma $-AG-groupoid $S,$ the following statements are true$.$
\end{theorem}

$(i)$ $A_{\alpha }$ is a $\Gamma $-right (left, two-sided) ideal of $S$ if $%
A $ is an intuitionistic fuzzy $\Gamma $-right (left) ideal of $S$.

$(ii)$ $A_{\alpha }$ is a $\Gamma $-bi-(generalized bi-) ideal of $S$ if $A$
is an intuitionistic fuzzy $\Gamma $-bi-(generalized bi-) ideal of $S$.

\begin{proof}
$(i)$: Let $S$ be a $\Gamma $-AG-groupoid and let $A$ be an intuitionistic
fuzzy $\Gamma $-right ideal of $S.$ If $x,y\in S\ $such that $x\in A_{\alpha
},$ then $\mu _{A}(x)\geq \alpha $ and $\gamma _{A}(x)\leq \alpha ,$
therefore $\mu _{A}(x\beta y)\geq \mu _{A}(x)\geq \alpha $ and $\gamma
_{A}(x\beta y)\leq \gamma _{A}(x)\leq \alpha ,$\ for all $\beta \in \Gamma $%
. Thus $x\beta y\in A_{\alpha },$ which shows that $A_{\alpha }$ is a $%
\Gamma $-right ideal of $S.$ Let $y\in A_{\alpha },$ then $\mu _{A}(y)\geq
\alpha $ and $\gamma _{A}(y)\leq \alpha .$ If $A$ is an intuitionistic fuzzy 
$\Gamma $-left ideal of $S,$ then $\mu _{A}(x\beta y)\geq \mu _{A}(y)\geq
\alpha $ and $\gamma _{A}(x\beta y)\leq \gamma _{A}(y)\leq \alpha $ implies
that $x\beta y\in A_{\alpha },$for all $\beta \in \Gamma .$ Which shows that 
$A_{\alpha }$ is a $\Gamma $-left ideal of $S.$

$(ii)$: Let $S$ be a $\Gamma $-AG-groupoid and let $A$ be an intuitionistic
fuzzy $\Gamma $-bi-(generalized bi-) ideal of $S.$ If $x,y$ and $z\in S$
such that $x$ and $z\in A_{\alpha },$ then $\mu _{A}(x)\geq \alpha ,$ $%
\gamma _{A}(x)\leq \alpha $, $\mu _{A}(z)\geq \alpha $ and $\gamma
_{A}(z)\leq \alpha $. Therefore, for all $\beta ,\delta \in \Gamma ,$ $\mu
_{A}((x\beta y)\delta z)\geq \mu _{A}(x)\wedge \mu _{A}(z)\geq \alpha ,$ and 
$\gamma _{A}((x\beta y)\delta z)\leq \gamma _{A}(x)\vee \gamma _{A}(z)\leq
\alpha $ implies that $(x\beta y)\delta z\in A_{\alpha },$ which shows that $%
A_{\alpha }$ is a generalized $\Gamma $-bi-ideal of $S$. Now let $x,y\in
A_{\alpha },$ then $\mu _{A}(x)\geq \alpha ,$ $\gamma _{A}(x)\leq \alpha ,$ $%
\mu _{A}(y)\geq \alpha $ and $\gamma _{A}(y)\leq \alpha $. Therefore for all 
$\beta \in \Gamma ,$ $\mu _{A}(x\beta y)\geq \mu _{A}(x)\wedge \mu
_{A}(y)\geq \alpha $ and $\gamma _{A}(x\beta y)\leq \gamma _{A}(x)\vee
\gamma _{A}(y)\leq \alpha ,$ implies that $x\beta y\in A_{\alpha }.$ Thus $%
A_{\alpha }$ is a $\Gamma $-bi-ideal of $S$.
\end{proof}

Conversely, let us define an $IFS$ $A=(\mu _{A},\gamma _{A})$ of an
AG-groupoid $S$ in Example \ref{ex} as follows: $\mu _{A}(1)=0.4,$ $\mu
_{A}(2)=0.8,$ $\mu _{A}(3)=\mu _{A}(4)=\mu _{A}(5)=0,$ $\gamma _{A}(1)=0.4,$ 
$\gamma _{A}(2)=0.3,$ $\gamma _{A}(3)=\gamma _{A}(4)=0.9$ and $\gamma
_{A}(5)=1.$ Let $\alpha =0.4,$ then it is easy to see that $A_{\alpha
}=\{a,b\}$ and one can easily verify from Example \ref{ex} that $\{a,b\}$ is
a right (left, bi-, generalized bi-) ideal of $S$ but $\mu _{A}(2\alpha
1)\ngeq \mu _{A}(2)$ $(\mu _{A}(1\alpha 2)\ngeq \mu _{A}(2),$ $\mu
_{A}((2\alpha 1)\alpha 2)\ngeq \mu _{A}(2))$ and $\gamma _{A}(2\alpha
1)\nleq \gamma _{A}(2)$ $(\gamma _{A}(1\alpha 2)\nleq \gamma _{A}(2),$ $%
\gamma _{A}((2\alpha 1)\alpha 2)\nleq \gamma _{A}(2))$ implies that $A$ is
not an intuitionistic fuzzy $\Gamma $-right (left, bi-, generalized bi-)
ideal of $S.$

Let $A=(\mu _{A},\gamma _{A})$ and $B=(\mu _{B},\gamma _{B})$ be any two $%
IFSs$ of a $\Gamma $-AG-groupoid $S$, then the product $A\circ _{\Gamma }B$
is defined by,

\begin{equation*}
\left( \mu _{A}\circ _{\Gamma }\mu _{B}\right) (a)=\left\{ 
\begin{array}{c}
\dbigvee\limits_{a=b\beta c}\left\{ \mu _{A}(b)\wedge \mu _{B}(c)\right\} 
\text{, if }a=b\beta c\text{ for some }b,\text{ }c\in S\text{ and }\beta \in
\Gamma . \\ 
0,\text{ otherwise.}%
\end{array}%
\right.
\end{equation*}

\begin{equation*}
\left( \gamma _{A}\circ _{\Gamma }\gamma _{B}\right) (a)=\left\{ 
\begin{array}{c}
\dbigwedge\limits_{a=b\beta c}\left\{ \gamma _{A}(b)\vee \gamma
_{B}(c)\right\} \text{, if }a=b\beta c\text{ for some }b,\text{ }c\in S\text{
and }\beta \in \Gamma . \\ 
1,\text{ otherwise.}%
\end{array}%
\right.
\end{equation*}

$A\subseteq B$ means that

\begin{center}
$\mu _{A}(x)\leq \mu _{B}(x)$ and $\gamma _{A}(x)\geq \gamma _{B}(x)$ for
all $x$ in $S.$
\end{center}

By keeping the generalization, the proof of the following is same as in \cite%
{Mordeson} and \cite{mad}.

\begin{lemma}
\label{as}Let $S$ be a $\Gamma $-AG-groupoid$,$ then the following holds.
\end{lemma}

$(i)$ An $IFS$ $A=(\mu _{A},\gamma _{A})$ is an intuitionistic fuzzy $\Gamma 
$-AG-subgroupoid of $S$ if and only if $\mu _{A}\circ _{\Gamma }\mu
_{A}\subseteq \mu _{A}$ and $\gamma _{A}\circ _{\Gamma }\gamma _{A}\supseteq
\gamma _{A}.$

$(ii)$ An $IFS$ $A=(\mu _{A},\gamma _{A})$ is intuitionistic fuzzy $\Gamma $%
-left (right) ideal of $S$ if and only if $S\circ _{\Gamma }\mu
_{A}\subseteq \mu _{A}$ and $\Theta \circ _{\Gamma }\gamma _{A}\supseteq
\gamma _{A}$ $(\mu _{A}\circ _{\Gamma }S\subseteq \mu _{A}$ and $\gamma
_{A}\circ _{\Gamma }\Theta \supseteq \gamma _{A}).$

\begin{theorem}
Let $A=(\mu _{A},\gamma _{A})$ be an $IFS$ of an intra-regular $\Gamma $-AG$%
^{\ast \ast }$-groupoid $S,$ then the following conditions are equivalent.
\end{theorem}

$(i)$ $A=(\mu _{A},\gamma _{A})$ is an intuitionistic fuzzy $\Gamma $%
-bi-ideal of $S$.

$(ii)$ $(A\circ _{\Gamma }\delta )\circ _{\Gamma }A=A$ and $A\circ _{\Gamma
}A=A,$ where $\delta =(S_{\delta },\Theta _{\delta }).$

\begin{proof}
$(i)\Longrightarrow (ii):$ Let $A=(\mu _{A},\gamma _{A})$ be an
intuitionistic fuzzy $\Gamma $-bi-ideal of an intra-regular $\Gamma $-AG$%
^{\ast \ast }$-groupoid $S$. Let $a\in A$, then there exists $x,$ $y\in S$
and $\alpha ,\beta ,\delta \in \Gamma $ such that $a=(x\alpha (a\beta
a))\delta y.$ Now let $\xi \in \Gamma ,$ then by using $(3),$ $(1),$ $(5),$ $%
(4)$ and $(2),$ we have%
\begin{eqnarray*}
a &=&(x\alpha (a\beta a))\delta y=(a\alpha (x\beta a))\delta y=(y\alpha
(x\beta a))\delta a=(y\alpha (x\beta ((x\alpha (a\beta a))\delta y)))\delta a
\\
&=&((u\xi v)\alpha (x\beta ((a\alpha (x\beta a))\delta y)))\delta
a=((((a\alpha (x\beta a))\delta y)\xi x)\alpha (v\delta u))\delta a \\
&=&(((x\delta y)\xi (a\alpha (x\beta a)))\alpha (v\delta u))\delta a=((a\xi
((x\delta y)\alpha (x\beta a)))\alpha (v\delta u))\delta a \\
&=&((a\xi ((x\delta x)\alpha (y\beta a)))\alpha (v\delta u))\delta
a=(((v\delta u)\xi ((x\delta x)\alpha (y\beta a)))\alpha a)\delta a \\
&=&(((v\delta u)\xi ((x\delta x)\alpha (y\beta ((x\alpha (a\beta a))\delta
y)))))\alpha a)\delta a \\
&=&(((v\delta u)\xi ((x\delta x)\alpha (y\beta ((x\alpha (a\beta a))\delta
(u\xi v)))))\alpha a)\delta a \\
&=&(((v\delta u)\xi ((x\delta x)\alpha (y\beta ((v\alpha u)\delta ((a\beta
a)\xi x)))))\alpha a)\delta a \\
&=&(((v\delta u)\xi ((x\delta x)\alpha (y\beta ((a\beta a)\delta ((v\alpha
u)\xi x)))))\alpha a)\delta a \\
&=&(((v\delta u)\xi ((x\delta x)\alpha ((a\beta a)\beta (y\delta ((v\alpha
u)\xi x)))))\alpha a)\delta a \\
&=&(((v\delta u)\xi ((a\beta a)\alpha ((x\delta x)\beta (y\delta ((v\alpha
u)\xi x)))))\alpha a)\delta a \\
&=&(((a\beta a)\xi ((v\delta u)\alpha ((x\delta x)\beta (y\delta ((v\alpha
u)\xi x)))))\alpha a)\delta a \\
&=&(((((x\delta x)\beta (y\delta ((v\alpha u)\xi x)))\beta (v\delta u))\xi
(a\alpha a))\alpha a)\delta a \\
&=&((a\xi ((((x\delta x)\beta (y\delta ((v\alpha u)\xi x)))\beta (v\delta
u))\alpha a))\alpha a)\delta a=(p\alpha a)\delta a
\end{eqnarray*}

where $p=a\xi ((((x\delta x)\beta (y\delta ((v\alpha u)\xi x)))\beta
(v\delta u))\alpha a).$ Therefore%
\begin{eqnarray*}
((\mu _{A}\circ _{\Gamma }S_{\delta })\circ _{\Gamma }\mu _{A})(a)
&=&\dbigvee\limits_{a=(p\alpha a)\delta a}\left\{ (\mu _{A}\circ _{\Gamma
}S_{\delta })(p\alpha a)\wedge \mu _{A}(a)\right\} \\
&\geq &\dbigvee\limits_{p\alpha a=p\alpha a}\left\{ \mu _{A}(p)\circ
_{\Gamma }S_{\delta }(a)\right\} \wedge \mu _{A}(a) \\
&\geq &\left\{ \mu _{A}(a\xi (((x^{2}\beta (y\delta ((v\alpha u)\xi
x)))\beta (v\delta u))\alpha a))\wedge S_{\delta }(a)\right\} \wedge \mu
_{A}(a) \\
&\geq &\mu _{A}(a)\wedge 1\wedge \mu _{A}(a)=\mu _{A}(a)
\end{eqnarray*}

and%
\begin{eqnarray*}
((\gamma _{A}\circ _{\Gamma }\Theta _{\delta })\circ _{\Gamma }\gamma
_{A})(a) &=&\dbigwedge\limits_{a=(p\alpha a)\delta a}\left\{ (\gamma
_{A}\circ _{\Gamma }\Theta _{\delta })(p\alpha a)\vee \gamma _{A}(a)\right\}
\\
&\leq &\dbigwedge\limits_{p\alpha a=p\alpha a}\left\{ \gamma _{A}(p)\circ
_{\Gamma }\Theta _{\delta }(a)\right\} \vee \gamma _{A}(a) \\
&\leq &\left\{ \gamma _{A}(a\xi (((x^{2}\beta (y\delta ((v\alpha u)\xi
x)))\beta (v\delta u))a))\vee \Theta _{\delta }(a)\right\} \vee \gamma
_{A}(a) \\
&\leq &\gamma _{A}(a)\vee 1\vee \gamma _{A}(a)=\gamma _{A}(a).
\end{eqnarray*}

This shows that $(\mu _{A}\circ _{\Gamma }S_{\delta })\circ _{\Gamma }\mu
_{A}\supseteq \mu _{A}$ and $(\gamma _{A}\circ _{\Gamma }\Theta _{\delta
})\circ _{\Gamma }\gamma _{A}\subseteq \gamma _{A},$ which implies that $%
(A\circ _{\Gamma }\delta )\circ _{\Gamma }A\supseteq A.$ Now by using $(4),$ 
$(1)$ and $(3),$ we have%
\begin{eqnarray*}
a &=&(x\alpha (a\beta a))\delta y=(a\alpha (x\beta a))\delta y=(y\alpha
(x\beta a))\delta a=(y\alpha (x\beta ((x\alpha (a\beta a))\delta y))\delta a
\\
&=&(y\alpha (x\beta ((x\alpha (a\beta a))\delta (u\xi v))))\delta a=(y\alpha
(x\beta ((v\alpha u)\delta ((a\beta a)\xi x))))\delta a \\
&=&(y\alpha (x\beta ((a\beta a)\delta ((v\alpha u)\xi x))))\delta a=(y\alpha
((a\beta a)\beta (x\delta ((v\alpha u)\xi x))))\delta a \\
&=&((a\beta a)\alpha (y\beta (x\delta ((v\alpha u)\xi x))))\delta
a=(((x\delta ((v\alpha u)\xi x))\beta y)\alpha (a\beta a))\delta a \\
&=&(a\alpha (((x\delta ((v\alpha u)\xi x))\beta y)\beta a))\delta a=(a\alpha
p)\delta a
\end{eqnarray*}

where $p=((x\delta ((v\alpha u)\xi x))\beta y)\beta a.$ Therefore%
\begin{eqnarray*}
((\mu _{A}\circ _{\Gamma }S_{\delta })\circ _{\Gamma }\mu _{A})(a)
&=&\dbigvee\limits_{a=(a\alpha p)\delta a}\left\{ (\mu _{A}\circ _{\Gamma
}S_{\delta })(a\alpha p)\wedge \mu _{A}(a)\right\} \\
&=&\dbigvee\limits_{a=(a\alpha p)\delta a}\left( \dbigvee\limits_{a\alpha
p=a\alpha p}\mu _{A}(a)\wedge S_{\delta }(p)\right) \wedge \mu _{A}(a) \\
&=&\dbigvee\limits_{a=(a\alpha p)\delta a}\left\{ \mu _{A}(a)\wedge 1\wedge
\mu _{A}(a)\right\} =\dbigvee\limits_{a=(ap)a}\mu _{A}(a)\wedge \mu _{A}(a)
\\
&\leq &\dbigvee\limits_{a=(a\alpha p)\delta a}\mu _{A}(((a\alpha ((x\delta
((v\alpha u)\xi x))\beta y)\beta a))\delta a)=\mu _{A}(a).
\end{eqnarray*}

This shows that $(\mu _{A}\circ _{\Gamma }S_{\delta })\circ _{\Gamma }\mu
_{A}\subseteq \mu _{A}$ and similarly we can show that $(\gamma _{A}\circ
_{\Gamma }\Theta _{\delta })\circ _{\Gamma }\gamma _{A}\supseteq \gamma
_{A}, $ which implies that $(A\circ _{\Gamma }\delta )\circ _{\Gamma
}A\subseteq A. $ Thus $(A\circ _{\Gamma }\delta )\circ _{\Gamma }A=A.$ We
have shown that $a=((a\xi ((((x\delta x)\beta (y\delta ((v\alpha u)\xi
x)))\beta (v\delta u))\alpha a))\alpha a)\delta a.$ Let $a=p\delta a$ where $%
p=(a\xi ((((x\delta x)\beta (y\delta ((v\alpha u)\xi x)))\beta (v\delta
u))\alpha a))\alpha a.$ Therefore 
\begin{eqnarray*}
(\mu _{A}\circ _{\Gamma }\mu _{A})(a) &=&\dbigvee\limits_{a=p\delta
a}\left\{ \mu _{A}(a\xi ((((x\delta x)\beta (y\delta ((v\alpha u)\xi
x)))\beta (v\delta u))\alpha a))\alpha a\wedge \mu _{A}(a)\right\} \\
&\geq &\mu _{A}(a)\wedge \mu _{A}(a)\wedge \mu _{A}(a)=\mu _{A}(a).
\end{eqnarray*}

This shows that $\mu _{A}\circ _{\Gamma }\mu _{A}\supseteq \mu _{A}$ and
similarly we can show that $\gamma _{A}\circ _{\Gamma }\gamma _{A}\subseteq
\gamma _{A}.$ Now by using Lemma \ref{as}, we get $A\circ _{\Gamma }A=A.$

$(ii)\Longrightarrow (i):$ Let $A=(\mu _{A},\gamma _{A})$ be an $IFS$ of an
intra-regular $\Gamma $-AG$^{\ast \ast }$-groupoid $S$, then let $\beta
,\delta \in \Gamma $, we have%
\begin{eqnarray*}
\mu _{A}((x\beta a)\delta y) &=&((\mu _{A}\circ _{\Gamma }S_{\delta })\circ
_{\Gamma }\mu _{A})((x\beta a)\delta y)=\dbigvee\limits_{(x\beta a)\delta
y=(x\beta a)\delta y}\{(\mu _{A}\circ _{\Gamma }S_{\delta })(x\beta a)\wedge
\mu _{A}(y)\} \\
&\geq &\dbigvee\limits_{x\beta a=x\beta a}\left\{ \mu _{A}(x)\wedge
S_{\delta }(a)\right\} \wedge \mu _{A}(y)\geq \mu _{A}(x)\wedge 1\wedge \mu
_{A}(y)=\mu _{A}(x)\wedge \mu _{A}(z).
\end{eqnarray*}

This shows that $\mu _{A}((x\beta a)\delta y)\geq \mu _{A}(x)\wedge \mu
_{A}(z)$ and similarly we can show that $\gamma _{A}((x\beta a)\delta y)\leq
\gamma _{A}(x)\vee \gamma _{A}(z)$. Also by Lemma \ref{as}, $A$ is an
intuitionistic fuzzy $\Gamma $-AG$^{\ast \ast }$-groupoid of $S$ and
therefore $A$ is an intuitionistic $\Gamma $-fuzzy bi-ideal of $S.$
\end{proof}

\begin{theorem}
Let $A=(\mu _{A},\gamma _{A})$ be an $IFS$ of an intra-regular $\Gamma $-AG$%
^{\ast \ast }$-groupoid $S,$ then the following conditions are equivalent.
\end{theorem}

$(i)$ $A=(\mu _{A},\gamma _{A})$ is an intuitionistic fuzzy $\Gamma $%
-interior ideal of $S$.

$(ii)$ $(\delta \circ _{\Gamma }A)\circ _{\Gamma }\delta =A,$ where $\delta
=(S_{\delta },\Theta _{\delta }).$

\begin{proof}
$(i)\Longrightarrow (ii):$ Let $A=(\mu _{A},\gamma _{A})$ be an
intuitionistic fuzzy $\Gamma $-interior ideal of an intra-regular $\Gamma $%
-AG$^{\ast \ast }$-groupoid . Let $a\in A$, then there exists $x,$ $y\in S$
and $\alpha ,\beta ,\delta \in \Gamma $, such that $a=(x\alpha (a\beta
a))\delta y.$ Now let $\xi \in \Gamma ,$ then by using $(3),$ $(5),$ $(4)$
and $(1),$ we have%
\begin{eqnarray*}
a &=&(x\alpha (a\beta a))\delta y=(a\beta (x\beta a))\delta y=((u\xi v)\beta
(x\beta a))\delta y \\
&=&((a\xi x)\beta (v\beta u))\delta y=(((v\beta u)\xi x)\beta a)\delta y.
\end{eqnarray*}

Therefore%
\begin{eqnarray*}
((S_{\delta }\circ _{\Gamma }\mu _{A})\circ _{\Gamma }S_{\delta })(a)
&=&\dbigvee\limits_{a=(((v\beta u)\xi x)\beta a)\delta y}\left\{ (S_{\delta
}\circ _{\Gamma }\mu _{A})(((v\beta u)\xi x)\beta a)\wedge S_{\delta
}(y)\right\} \\
&\geq &\dbigvee\limits_{((v\beta u)\xi x)\beta a=((v\beta u)\xi x)\beta
a}\left\{ (S_{\delta }((v\beta u)\xi x)\wedge \mu _{A}(a)\right\} \wedge 1 \\
&\geq &1\wedge \mu _{A}(a)\wedge 1=\mu _{A}(a).
\end{eqnarray*}

This proves that $(S_{\delta }\circ _{\Gamma }\mu _{A})\circ _{\Gamma
}S_{\delta }\supseteq \mu _{A}$ and similarly we can show that $(\Theta
_{\delta }\circ _{\Gamma }\gamma _{A})\circ _{\Gamma }\Theta _{\delta
}\subseteq \gamma _{A}$, therefore $(\delta \circ _{\Gamma }A)\circ _{\Gamma
}\delta \supseteq A.$ Now again%
\begin{eqnarray*}
((S_{\delta }\circ _{\Gamma }\mu _{A})\circ _{\Gamma }S_{\delta })(a)
&=&\dbigvee\limits_{a=(x\alpha (a\beta a))\delta y}\left\{ (S_{\delta }\circ
_{\Gamma }\mu _{A})(x\alpha (a\beta a))\wedge S_{\delta }(y)\right\} \\
&=&\dbigvee\limits_{a=(x\alpha (a\beta a))\delta y}\left(
\dbigvee\limits_{x\alpha (a\beta a)=x\alpha (a\beta a)}S_{\delta }(x)\wedge
\mu _{A}(a\beta a)\right) \wedge S_{\delta }(y) \\
&=&\dbigvee\limits_{a=(x\alpha (a\beta a))\delta y}\left\{ 1\wedge \mu
_{A}(a\beta a)\wedge 1\right\} =\dbigvee\limits_{a=(x\alpha (a\beta
a))\delta y}\mu _{A}(a\beta a) \\
&\leq &\dbigvee\limits_{a=(x\alpha (a\beta a))\delta y}\mu _{A}((x\alpha
(a\beta a))\delta y)=\mu _{A}(a).
\end{eqnarray*}

Thus $(S_{\delta }\circ _{\Gamma }\mu _{A})\circ _{\Gamma }S_{\delta
}\subseteq \mu _{A}$ and similarly we can show that $(\Theta _{\delta }\circ
_{\Gamma }\gamma _{A})\circ _{\Gamma }\Theta _{\delta }\supseteq \gamma
_{A}, $ therefore $(\delta \circ _{\Gamma }A)\circ _{\Gamma }\delta
\subseteq A.$ Hence it follows that $(\delta \circ _{\Gamma }A)\circ
_{\Gamma }\delta =A.$

$(ii)\Longrightarrow (i):$ Let $A=(\mu _{A},\gamma _{A})$ be an $IFS$ of an
intra-regular $\Gamma $-AG$^{\ast \ast }$-groupoid $S$, and let $\beta
,\delta \in \Gamma ,$ then%
\begin{eqnarray*}
\mu _{A}((x\beta a)\delta y) &=&((S_{\delta }\circ _{\Gamma }\mu _{A})\circ
_{\Gamma }S_{\delta })((x\beta a)\delta y)=\dbigvee\limits_{(x\beta a)\delta
y=(x\beta a)\delta y}\left\{ (S_{\delta }\circ _{\Gamma }\mu _{A})(x\beta
a)\wedge S_{\delta }(y)\right\} \\
&\geq &\dbigvee\limits_{x\beta a=x\beta a}\left\{ (S_{\delta }(x)\circ
_{\Gamma }\mu _{A}(a)\right\} \wedge S_{\delta }(y)\geq \mu _{A}(a).
\end{eqnarray*}

Similarly we can show that $\gamma _{A}((x\beta a)\delta y)\leq \gamma
_{A}(a)$ and therefore $A=(\mu _{A},\gamma _{A})$ is an intuitionistic fuzzy 
$\Gamma $-interior ideal of $S.$
\end{proof}

\begin{lemma}
\label{iff}Let $A=(\mu _{A},\gamma _{A})$ be an $IFS$ of an intra-regular $%
\Gamma $-AG$^{\ast \ast }$-groupoid$,$ then $A=(\mu _{A},\gamma _{A})$ is an
intuitionistic fuzzy $\Gamma $-left ideal of $S$ if and only if $A=(\mu
_{A},\gamma _{A})$ is an intuitionistic fuzzy $\Gamma $-right ideal of $S.$
\end{lemma}

\begin{proof}
Let $S$ be an intra-regular $\Gamma $-AG$^{\ast \ast }$-groupoid and let $%
A=(\mu _{A},\gamma _{A})$ be an intuitionistic fuzzy $\Gamma $-left ideal of 
$S.$ Now for $a,b\in S$ there exists $x,y,x^{^{\prime }},y^{^{\prime }}\in S$
and $\alpha ,\beta ,\delta \in \Gamma $ such that $a=(x\alpha (a\beta
a))\delta y$ and $b=(x^{^{\prime }}b^{2})y^{^{\prime }}.$ Now let $\xi ,\psi
\in \Gamma ,$ then by using $(1),$ $(4),$ $(5)$ and $(3),$ we have 
\begin{eqnarray*}
\mu _{A}(a\xi b) &=&\mu _{A}(((x\alpha (a\beta a))\delta y)\xi b)=\mu
_{A}((b\delta y)\xi (x\alpha (a\beta a)))=\mu _{A}(((a\beta a)\delta x)\xi
(y\alpha b)) \\
&=&\mu _{A}(((x\beta a)\delta a)\xi (y\alpha b))=\mu _{A}(((x\beta a)\delta
(u\psi v))\xi (y\alpha b)) \\
&=&\mu _{A}(((v\beta u)\delta (a\psi x))\xi (y\alpha b))=\mu _{A}((a\delta
((v\beta u)\psi x))\xi (y\alpha b)) \\
&=&\mu _{A}(((y\alpha b)\delta ((v\beta u)\psi x))\xi a)\geq \mu _{A}(a).
\end{eqnarray*}

Similarly we can get $\gamma _{A}(a\xi b)\leq \gamma _{A}(a),$ which implies
that $A=(\mu _{A},\gamma _{A})$ is an intuitionistic fuzzy $\Gamma $-right
ideal of $S.$

Conversely let $A=(\mu _{A},\gamma _{A})$ be an intuitionistic fuzzy $\Gamma 
$-right ideal of $S.$ Now by using $(3)$ and $(4),$ we have 
\begin{eqnarray*}
\mu _{A}(ab) &=&\mu _{A}(a((x^{^{\prime }}b^{2})y^{^{\prime }})=\mu
_{A}((x^{^{\prime }}b^{2})(ay^{^{\prime }}))=\mu _{A}((y^{^{\prime
}}a)(b^{2}x^{^{\prime }})) \\
&=&\mu _{A}(b^{2}((y^{^{\prime }}a)x))\geq \mu _{A}(b).
\end{eqnarray*}

Also we can get $\gamma _{A}(ab)\leq \gamma _{A}(b),$ which implies that $%
A=(\mu _{A},\gamma _{A})$ is an intuitionistic fuzzy $\Gamma $-left ideal of 
$S.$
\end{proof}

A $\Gamma $-AG-groupoid $S$ is called a $\Gamma $-left (right) duo if every $%
\Gamma $-left (right) ideal of $S$ is a $\Gamma $-two-sided ideal of $S$ and
is called a $\Gamma $-duo if it is both a $\Gamma $-left and a $\Gamma $%
-right duo.

A $\Gamma $-AG-groupoid $S$ is called an intuitionistic fuzzy $\Gamma $-left
(right) duo if every intuitionistic fuzzy $\Gamma $-left (right) ideal of $S$
is an intuitionistic fuzzy $\Gamma $-two-sided ideal of $S$ and is called an
intuitionistic fuzzy $\Gamma $-duo if it is both an intuitionistic fuzzy $%
\Gamma $-left and an intuitionistic fuzzy $\Gamma $-right duo.

\begin{corollary}
Every intra-regular $\Gamma $-AG$^{\ast \ast }$-groupoid $S$ is an
intuitionistic fuzzy $\Gamma $-duo.
\end{corollary}

\begin{lemma}
\label{aqw}In an intra-regular $\Gamma $-AG-groupoid $S,$ $\delta \circ
_{\Gamma }A=A$ and $A\circ _{\Gamma }\delta =A$ holds for an $IFS$ $A=(\mu
_{A},\gamma _{A})$ of $S$ where $\delta =(S_{\delta },\Theta _{\delta }).$
\end{lemma}

\begin{proof}
Let $A=(\mu _{A},\gamma _{A})$ be an $IFS$ of an intra-regular $\Gamma $%
-AG-groupoid $S$ and let $a\in S,$ then there exist $x,y\in S$ and $\alpha
,\beta ,\delta \in \Gamma $ such that $a=(x\alpha (a\beta a))\gamma y.$ Now
by using $(3)$ and $(1),$ we have%
\begin{equation*}
a=(x\alpha (a\beta a))\delta y=(a\alpha (x\beta a))\delta y=(y\alpha (x\beta
a))\delta a.
\end{equation*}

Therefore%
\begin{eqnarray*}
(S_{\delta }\circ _{\Gamma }\mu _{A})(a) &=&\dbigvee\limits_{a=(y\alpha
(x\beta a))\delta a}\{S_{\delta }(y\alpha (x\beta a))\wedge \mu
_{A}(a)\}=\dbigvee\limits_{a=(y\alpha (x\beta a))\delta a}\{1\wedge \mu
_{A}(a)\} \\
&=&\dbigvee\limits_{a=(y\alpha (x\beta a))\delta a}\mu _{A}(a)=\mu _{A}(a).
\end{eqnarray*}

Similarly we can show that $\Theta _{\delta }\circ _{\Gamma }\gamma
_{A}=\gamma _{A},$ which shows that $\delta \circ _{\Gamma }A=A.$

Now let $\xi \in \Gamma ,$ then by using $(5),$ $(4)$ and $(3),$ we have%
\begin{eqnarray*}
a &=&(x\alpha (a\beta a))\delta y=(x\alpha (a\beta a))\delta (u\xi
v)=(v\alpha u)\delta ((a\beta a)\xi x) \\
&=&(a\beta a)\delta ((v\alpha u)\xi x)=(x\beta (v\alpha u))\delta (a\xi
a)=a\delta ((x\beta (v\alpha u))\xi a).
\end{eqnarray*}

Therefore%
\begin{eqnarray*}
(\mu _{A}\circ _{\Gamma }S_{\delta })(a) &=&\dbigvee\limits_{a=a\delta
((x\beta (v\alpha u))\xi a)}\{\mu _{A}(a)\wedge S_{\delta }((x\beta (v\alpha
u))\xi a)\}=\dbigvee\limits_{a=a\delta ((x\beta (v\alpha u))\xi a)}\{\mu
_{A}(a)\wedge 1\} \\
&=&\dbigvee\limits_{a=a\delta ((x\beta (v\alpha u))\xi a)}\mu _{A}(a)=\mu
_{A}(a).
\end{eqnarray*}

Similarly we can show that $\gamma _{A}\circ _{\Gamma }\Theta _{\delta
}=\gamma _{A}$ which shows that $A\circ _{\Gamma }\delta =A.$
\end{proof}

\begin{corollary}
\label{cor}In an intra-regular $\Gamma $-AG-groupoid $S,$ $\delta \circ
_{\Gamma }A=A$ and $A\circ _{\Gamma }\delta =A$ holds for every
intuitionistic fuzzy $\Gamma $-left $($right$,$ two-sided$)$ $A=(\mu
_{A},\gamma _{A})$ of $S,$ where $\delta =(S_{\delta },\Theta _{\delta }).$
\end{corollary}

\begin{lemma}
\label{ss}In an intra-regular $\Gamma $-AG-groupoid $S,$ $\delta \circ
_{\Gamma }\delta =\delta ,$ where $\delta =(S_{\delta },\Theta _{\delta }).$
\end{lemma}

\begin{proof}
Let $S$ be an intra-regular $\Gamma $-AG-groupoid, then%
\begin{equation*}
(S_{\delta }\circ _{\Gamma }S_{\delta })(a)=\dbigvee\limits_{a=(x\alpha
(a\beta a))\delta y}\{S_{\delta }(x\alpha (a\beta a))\wedge S_{\delta
}(y)\}=1=S_{\delta }(a)
\end{equation*}

and%
\begin{equation*}
(\Theta _{\delta }\circ _{\Gamma }\Theta _{\delta
})(a)=\dbigwedge\limits_{a=(x\alpha (a\beta a))\delta y}\{\Theta _{\delta
}(x\alpha (a\beta a))\vee \Theta _{\delta }(y)\}=0=\Theta _{\delta }(a).
\end{equation*}
\end{proof}

\begin{theorem}
\label{ae}Let $A=(\mu _{A},\gamma _{A})$ be an $IFS$ of an intra-regular $%
\Gamma $-AG$^{\ast \ast }$-groupoid $S,$ then the following conditions are
equivalent.
\end{theorem}

$(i)$ $A=(\mu _{A},\gamma _{A})$ is an intuitionistic fuzzy $\Gamma $-quasi
ideal of $S$.

$(ii)$ $(A\circ _{\Gamma }\delta )\cap (\delta \circ _{\Gamma }A)=A,$ where $%
\delta =(S_{\delta },\Theta _{\delta }).$

\begin{proof}
$(i)\Longrightarrow (ii)$ can be followed from Lemma \ref{aqw} and $%
(ii)\Longrightarrow (i)$ is obvious.
\end{proof}

\begin{theorem}
\label{qu}Let $A=(\mu _{A},\gamma _{A})$ an $IFS$ of an intra-regular $%
\Gamma $-AG$^{\ast \ast }$-groupoid $S$, then the following statements are
equivalent.
\end{theorem}

$(i)$ $A$ is an intuitionistic fuzzy $\Gamma $-two-sided ideal of $S.$

$(ii)$ $A$ is an intuitionistic fuzzy $\Gamma $-quasi ideal of $S.$

\begin{proof}
$(i)\Longrightarrow (ii)$ is an easy consequence of Corollary \ref{cor} and
Theorem \ref{ae}.

$(ii)\Longrightarrow (i):$ Let $A=(\mu _{A},\gamma _{A})$ be an
intuitionistic fuzzy $\Gamma $-quasi ideal of an intra-regular $\Gamma $-AG$%
^{\ast \ast }$-groupoid $S$ and let $a\in S,$ then there exist $x,y\in S$
and $\alpha ,\beta ,\delta \in \Gamma $ such that $a=(x\alpha (a\beta
a))\delta y.$ Now let $\xi ,\psi \in \Gamma ,$ then by using $(3),$ $(5)$
and $(4),$ we have%
\begin{eqnarray*}
a &=&(x\alpha (a\beta a))\delta y=(a\alpha (x\beta a))\delta (u\xi
v)=(v\alpha u)\delta ((x\beta a)\xi a) \\
&=&(v\alpha u)\delta ((x\beta a)\xi (m\psi n))=(v\alpha u)\delta ((n\beta
m)\xi (a\psi x)) \\
&=&(v\alpha u)\delta (a\xi ((n\beta m)\psi x))=a\delta ((v\alpha u)\xi
((n\beta m)\psi x)).
\end{eqnarray*}

Therefore%
\begin{eqnarray*}
(\mu _{A}\circ _{\Gamma }S_{\delta })(a) &=&\dbigvee\limits_{a=a\delta
((v\alpha u)\xi ((n\beta m)\psi x))}\{\mu _{A}(a)\wedge S_{\delta }((v\alpha
u)\xi ((n\beta m)\psi x))\} \\
&\geq &\mu _{A}(a)\wedge 1=\mu _{A}(a)\text{.}
\end{eqnarray*}%
Similarly we can show that $\gamma _{A}\circ _{\Gamma }\Theta _{\delta
}\subseteq \gamma _{A}$ which implies that $A\circ _{\Gamma }\delta
\supseteq A.$ Now by using Lemmas \ref{aqw}, \ref{ss} and $(2),$ we have%
\begin{equation*}
A\circ _{\Gamma }\delta =(\delta \circ _{\Gamma }A)\circ _{\Gamma }(\delta
\circ _{\Gamma }\delta )=(\delta \circ _{\Gamma }\delta )\circ _{\Gamma
}(A\circ _{\Gamma }\delta )=\delta \circ _{\Gamma }(A\circ _{\Gamma }\delta
)\supseteq \delta \circ _{\Gamma }A\text{.}
\end{equation*}

This shows that $\delta \circ _{\Gamma }A\subseteq (A\circ _{\Gamma }\delta
)\cap (\delta \circ _{\Gamma }A).$ As $A$ is an intuitionistic fuzzy $\Gamma 
$-quasi ideal of $S,$ thus we get $\delta \circ _{\Gamma }A\subseteq A$. Now
by using Lemma \ref{as}, $A$ is an intuitionistic fuzzy $\Gamma $-left ideal
of $S$ and by Lemma \ref{iff}, $A$ is an intuitionistic fuzzy $\Gamma $%
-right ideal of $S,$ that is, $A\ $is an intuitionistic fuzzy $\Gamma $%
-two-sided ideal of $S.$
\end{proof}

\begin{theorem}
\label{intr}Let $A=(\mu _{A},\gamma _{A})$ be an $IFS$ of an intra-regular $%
\Gamma $-AG$^{\ast \ast }$-groupoid $S$, then the following statements are
equivalent.
\end{theorem}

$(i)$ $A$ is an intuitionistic fuzzy $\Gamma $-two-sided ideal of $S.$

$(ii)$ $A$ is an intuitionistic fuzzy $\Gamma $-interior ideal of $S.$

\begin{proof}
$(i)\Longrightarrow (ii)$ is obvious.

$(ii)\Longrightarrow (i):$ Let $A=(\mu _{A},\gamma _{A})$ be an
intuitionistic fuzzy $\Gamma $-interior ideal of an intra-regular $\Gamma $%
-AG$^{\ast \ast }$-groupoid $S$ and let $a,b\in S,$ then there exist $x,y\in
S$ and $\alpha ,\beta ,\delta \in \Gamma $ such that $a=(x\alpha (a\beta
a))\delta y.$ Now let $\xi \in \Gamma ,$ then by using $(3)$, $(1)$ and $%
(4), $ we have%
\begin{eqnarray*}
\mu _{A}(a\xi b) &=&\mu _{A}(((x\alpha (a\beta a))\delta y)\xi b)=\mu
_{A}(((a\alpha (x\beta a))\delta y)\xi b)=\mu _{A}((b\delta y)\xi (a\alpha
(x\beta a))) \\
&=&\mu _{A}(((x\beta a)\delta a)\xi (y\alpha b))\geq \mu _{A}(a).
\end{eqnarray*}

Similarly we can prove that $\gamma _{A}(a\xi b)\leq \gamma _{A}(a).$ Thus $%
A $ is an intuitionistic fuzzy $\Gamma $-right ideal of $S$ and by using
Lemma \ref{iff}, $A$ is an intuitionistic fuzzy $\Gamma $-two-sided ideal of 
$S.$
\end{proof}

\begin{theorem}
Let $A=(\mu _{A},\gamma _{A})$ be an $IFS$ of an intra-regular $\Gamma $-AG$%
^{\ast \ast }$-groupoid $S$, then the following statements are equivalent.
\end{theorem}

$(i)$ $A$ is an intuitionistic fuzzy $\Gamma $-left ideal of $S$.

$(ii)$ $A$ is an intuitionistic fuzzy $\Gamma $-right ideal of $S$.

$(iii)$ $A$ is an intuitionistic fuzzy $\Gamma $-two-sided ideal of $S$.

$(iv)$ $A$ is an intuitionistic fuzzy $\Gamma $-bi-ideal of $S$.

$(v)$ $A$ is an intuitionistic fuzzy $\Gamma $-generalized bi-ideal of $S$.

$(vi)$ $A$ is an intuitionistic fuzzy $\Gamma $-interior ideal of $S$.

$(vii)$ $A$ is an intuitionistic fuzzy $\Gamma $-quasi ideal of $S.$

$(viii)$ $A\circ _{\Gamma }\delta =A$ and $\delta \circ _{\Gamma }A=A.$

\begin{proof}
$(i)\Longrightarrow (viii)$ can be followed from Corollary \ref{cor} and $%
(ix)\Longrightarrow (viii)$ is obvious.

$(vii)\Longrightarrow (vi):$ Let $A=(\mu _{A},\gamma _{A})$ be an
intuitionistic fuzzy $\Gamma $-quasi ideal of an intra-regular $\Gamma $-AG$%
^{\ast \ast }$-groupoid $S$. Now for $a\in S$ there exist $x,y\in S$ and $%
\alpha ,\beta ,\delta \in \Gamma $ such that $a=(b\alpha (a\beta a))\delta
c. $ Now let $\xi ,\psi \in \Gamma ,$ then by using $(3),$ $(4)$ and $(1)$,
we have%
\begin{eqnarray*}
(x\xi a)\psi y &=&(x\xi ((b\alpha (a\beta a))\delta c))\psi y=((b\alpha
(a\beta a))\xi (x\delta c))\psi y \\
&=&((c\alpha x)\xi ((a\beta a)\delta b))\psi y=((a\beta a)\xi ((c\alpha
x)\delta b))\psi y \\
&=&(y\xi ((c\alpha x)\delta b))\psi (a\beta a)=a\psi ((y\xi ((c\alpha
x)\delta b))\beta a)
\end{eqnarray*}

and%
\begin{eqnarray*}
(x\xi a)\psi y &=&(x\xi ((b\alpha (a\beta a))\delta c))\psi y=((b\alpha
(a\beta a))\xi (x\delta c))\psi y \\
&=&((c\alpha x)\xi ((a\beta a)\delta b))\psi y=((a\beta a)\xi ((c\alpha
x)\delta b))\psi y \\
&=&(y\xi ((c\alpha x)\delta b))\psi (a\beta a)=(a\xi a)\psi (((c\alpha
x)\delta b)\beta y) \\
&=&((((c\alpha x)\delta b)\beta y)\xi a)\psi a.
\end{eqnarray*}%
Now by using Theorem \ref{ae}, we have%
\begin{equation*}
\mu _{A}((x\xi a)\psi y)=((\mu _{A}\circ _{\Gamma }S_{\delta })\cap
(S_{\delta }\circ _{\Gamma }\mu _{A}))((x\xi a)\psi y)=(\mu _{A}\circ
_{\Gamma }S_{\delta })((x\xi a)\psi y)\wedge (S_{\delta }\circ _{\Gamma }\mu
_{A})((x\xi a)\psi y).
\end{equation*}

Now%
\begin{equation*}
(\mu _{A}\circ _{\Gamma }S_{\delta })((x\xi a)\psi y)=\dbigvee\limits_{a\psi
((y\xi ((c\alpha x)\delta b))\beta a)}\left\{ \mu _{A}(a)\wedge S_{\delta
}((y\xi ((c\alpha x)\delta b))\beta a)\right\} \geq \mu _{A}(a)
\end{equation*}

and 
\begin{equation*}
\left( S_{\delta }\circ _{\Gamma }\mu _{A}\right) ((x\xi a)\psi
y)=\dbigvee\limits_{(xa)y=((((c\alpha x)\delta b)\beta y)\xi a)\psi
a}\left\{ S_{\delta }((((c\alpha x)\delta b)\beta y)\xi a)\wedge \mu
_{A}(a)\right\} \geq \mu _{A}(a).
\end{equation*}

This implies that $\mu _{A}((x\xi a)\psi y)\geq \mu _{A}(a)$ and similarly
we can show that $\gamma _{A}((x\xi a)\psi y)\leq \gamma _{A}(a)$. Thus $A$
is an intuitionistic fuzzy $\Gamma $-interior ideal of $S.$

$(vi)\Longrightarrow (v):$ Let $A$ be an intuitionistic fuzzy $\Gamma $%
-interior ideal of $S,$ then by Theorem \ref{intr}, $A$ is an intuitionistic
fuzzy $\Gamma $-two-sided ideal of $S$ and it is easy to observe that $A$ is
an intuitionistic fuzzy $\Gamma $-generalized bi-ideal of $S$.

$(v)\Longrightarrow (iv)$: Let $A=(\mu _{A},\gamma _{A})$ be an
intuitionistic fuzzy $\Gamma $-generalized bi-ideal of an intra-regular $%
\Gamma $-AG$^{\ast \ast }$-groupoid $S$ Let $a\in S$, then there exists $%
x,y\in S$ and $\alpha ,\beta ,\delta \in \Gamma $ such that $a=(x\alpha
(a\beta a))\delta y.$ Now let $\xi ,\psi ,\eta \in \Gamma ,$ then by using $%
(3),$ $(5),$ $(4)$ and $(1),$ we have 
\begin{eqnarray*}
\mu _{A}(a\xi b) &=&\mu _{A}(((x\alpha (a\beta a))\delta y)\xi b)=\mu
_{A}(((a\alpha (x\beta a))\delta y)\xi b) \\
&=&\mu _{A}((((s\psi t)\alpha (x\beta a))\delta y)\xi b)=\mu _{A}((((a\psi
x)\alpha (t\beta s))\delta y)\xi b)=\mu _{A}(((t\alpha ((a\psi x)\beta
s))\delta y)\xi b) \\
&=&\mu _{A}(((t\alpha ((a\psi x)\beta s))\delta (u\eta v))\xi b)=\mu
_{A}(((v\alpha u)\delta (((a\psi x)\beta s)\eta t))\xi b) \\
&=&\mu _{A}(((v\alpha u)\delta ((t\beta s)\eta (a\psi x)))\xi b)=\mu
_{A}(((v\alpha u)\delta (a\eta ((t\beta s)\psi x)))\xi b) \\
&=&\mu _{A}((a\delta ((v\alpha u)\eta ((t\beta s)\psi x)))\xi b)\geq \mu
_{A}(a)\wedge \mu _{A}(b).
\end{eqnarray*}

Similarly we can show that $\gamma _{A}(a\xi b)\leq \gamma _{A}(a)\vee
\gamma _{A}(b)$ and therefore $A=(\mu _{A},\gamma _{A})$ is an
intuitionistic fuzzy $\Gamma $-bi-ideal of $S$.

$(iv)\Longrightarrow (iii):$ Let $A=(\mu _{A},\gamma _{A})$ be an
intuitionistic fuzzy $\Gamma $-bi-ideal of an intra-regular $\Gamma $-AG$%
^{\ast \ast }$-groupoid $S$. Let $a,t\in S$, then there exists $%
x,y,x^{^{\prime }},y^{^{\prime }}\in S$ $\ $and $\alpha ,\beta ,\delta \in
\Gamma $ such that $a=(x\alpha (a\beta a))\delta y$ and $t=(x^{^{\prime
}}\alpha (a\beta a))\delta y^{^{\prime }}.$ Now let $\xi ,\gamma \in \Gamma
, $ then by using $(3),$ $(1)$, $(5)$ and $(4),$ we have%
\begin{eqnarray*}
\mu _{A}(a\xi b) &=&\mu _{A}(((x\alpha (a\beta a))\delta y)\xi b)=\mu
_{A}(((a\alpha (x\beta a))\delta y)\xi b)=\mu _{A}((b\delta y)\xi (a\alpha
(x\beta a))) \\
&=&\mu _{A}((b\delta y)\xi ((s\gamma t)\alpha (x\beta a)))=\mu _{A}((b\delta
y)\xi ((a\gamma x)\alpha (t\beta s))) \\
&=&\mu _{A}(((t\beta s)\delta (a\gamma x))\xi (y\alpha b))=\mu _{A}((a\delta
((t\beta s)\gamma x))\xi (y\alpha b)) \\
&=&\mu _{A}(((y\alpha b)\delta ((t\beta s)\gamma x))\xi a)=\mu
_{A}(((y\alpha b)\delta ((((x\alpha a^{2})\delta y)\beta s)\gamma x))\xi a)
\\
&=&\mu _{A}(((y\alpha b)\delta ((x\beta s)\gamma ((x\alpha a^{2})\delta
y)))\xi a)=\mu _{A}(((y\alpha b)\delta ((y\beta (x\alpha a^{2}))\gamma
(s\delta x)))\xi a) \\
&=&\mu _{A}(((y\alpha b)\delta ((x\beta s)\gamma ((x\alpha a^{2})\delta
y)))\xi a)=\mu _{A}(((yb)((xs)((x\alpha a^{2})\delta (uv))))\xi a) \\
&=&\mu _{A}(((yb)((xs)((vu)(a^{2}x))))\xi a)=\mu
_{A}(((yb)((xs)(a^{2}((vu)x))))\xi a) \\
&=&\mu _{A}(((yb)(a^{2}((xs)((vu)x))))\xi a)=\mu
_{A}((a^{2}((yb)((xs)((vu)x))))\xi a) \\
&\geq &\mu _{A}(a^{2})\wedge \mu _{A}(a)\geq \mu _{A}(a)\wedge \mu
_{A}(a)\wedge \mu _{A}(a)=\mu _{A}(a).
\end{eqnarray*}

Similarly we can prove that $\gamma _{A}(a\xi b)\leq \gamma _{A}(a)$ and
therefore $A=(\mu _{A},\gamma _{A})$ is an intuitionistic fuzzy $\Gamma $%
-right ideal of $S$. Now by using Lemma \ref{iff}, $A=(\mu _{A},\gamma _{A})$
is an intuitionistic fuzzy $\Gamma $-two-sided ideal of $S$.

$(iii)\Longrightarrow (ii)$ and $(ii)\Longrightarrow (i)$ are an easy
consequences of Lemma \ref{iff}.
\end{proof}

Let $A=(\mu _{A},\gamma _{A})$ and $B=(\mu _{B},\gamma _{B})$ are $IFSs$ of
an AG-groupoid $S.$ The symbols $A\cap B$ will means the following $IFS$ of $%
S$

\begin{center}
$(\mu _{A}\cap \mu _{B})(x)=\min \{\mu _{A}(x),\mu _{B}(x)\}=\mu
_{A}(x)\wedge \mu _{B}(x),$ for all $x$ in $S.$

$(\gamma _{A}\cup \gamma _{B})(x)=\max \{\gamma _{A}(x),\gamma
_{B}(x)\}=\gamma _{A}(x)\vee \gamma _{B}(x),$ for all $x$ in $S.$
\end{center}

The symbols $A\cup B$ will means the following $IFS$ of $S$

\begin{center}
$(\mu _{A}\cup \mu _{B})(x)=\max \{\mu _{A}(x),\mu _{B}(x)\}=\mu _{A}(x)\vee
\mu _{B}(x),$ for all $x$ in $S.$

$(\gamma _{A}\cap \gamma _{B})(x)=\min \{\gamma _{A}(x),\gamma
_{B}(x)\}=\gamma _{A}(x)\wedge \gamma _{B}(x),$ for all $x$ in $S.$
\end{center}

\begin{lemma}
\label{fgh}Let $S$ be an intra-regular $\Gamma $-AG$^{\ast \ast }$-groupoid
and let $A=(\mu _{A},\gamma _{A})$ and $B=(\mu _{B},\gamma _{B})$ are any
intuitionistic fuzzy $\Gamma $-two-sided ideals of $S,$ then $A\circ
_{\Gamma }B=A\cap B$.
\end{lemma}

\begin{proof}
Assume that $A=(\mu _{A},\gamma _{A})$ and $B=(\mu _{B},\gamma _{B})$ are
any intuitionistic fuzzy $\Gamma $-two-sided ideals of an intra-regular AG$%
^{\ast \ast }$-groupoid $S$, then by using Lemma \ref{as}, we have $\mu
_{A}\circ _{\Gamma }\mu _{B}\subseteq \mu _{A}\cap \mu _{B}$ and $\gamma
_{A}\circ _{\Gamma }\gamma _{B}\supseteq \gamma _{A}\cup \gamma _{B},$ which
shows that $A\circ _{\Gamma }B\subseteq A\cap B$. Let $a\in S,$ then there
exists $x,y\in S$ and $\alpha ,\beta ,\delta \in \Gamma $ such that $%
a=(x\alpha (a\beta a))\delta y.$ Now let $\xi \in \Gamma ,$ then by using $%
(3),$ $(5)$ and $(2),$ we have%
\begin{equation*}
a=(x\alpha (a\beta a))\delta y=(a\alpha (x\beta a))\delta (u\xi v)=(a\alpha
u)\delta ((x\beta a)\xi v).
\end{equation*}%
Therefore, we have%
\begin{eqnarray*}
(\mu _{A}\circ _{\Gamma }\mu _{B})(a) &=&\dbigvee\limits_{a=(a\alpha
u)\delta ((x\beta a)\xi v)}\{\mu _{A}(a\alpha u)\wedge \mu _{B}((x\beta
a)\xi v)\}\geq \mu _{A}(a\alpha u)\wedge \mu _{B}((x\alpha a)\xi v) \\
&\geq &\mu _{A}(a)\wedge \mu _{B}(a)=(\mu _{A}\cap \mu _{B})(a)
\end{eqnarray*}

and%
\begin{eqnarray*}
(\gamma _{A}\circ _{\Gamma }\gamma _{A})(a) &=&\dbigwedge\limits_{a=(a\alpha
u)\delta ((x\beta a)\xi v)}\left\{ \gamma _{A}(a\alpha u)\vee \gamma
_{A}((x\beta a)\xi v)\right\} \leq \gamma _{A}(a\alpha u)\vee \gamma
_{A}((x\alpha a)\xi v) \\
&\leq &\gamma _{A}(a)\vee \gamma _{A}(a)=(\gamma _{A}\cup \gamma _{A})(a).
\end{eqnarray*}

Thus we get that $\mu _{A}\circ _{\Gamma }\mu _{B}\supseteq \mu _{A}\cap \mu
_{B}$ and $\gamma _{A}\circ _{\Gamma }\gamma _{B}\subseteq \gamma _{A}\cup
\gamma _{B},$ which give us $A\circ _{\Gamma }B\supseteq A\cap B$ and
therefore $A\circ _{\Gamma }B=A\cap B.$
\end{proof}

The converse of Lemma \ref{fgh} is not true in general which is discussed in
the following.

Let us consider an AG-groupoid $S=\left\{ 1,2,3,4,5\right\} $ in the
following Cayley's table.

\begin{center}
\begin{tabular}{l|lllll}
. & $1$ & $2$ & $3$ & $4$ & $5$ \\ \hline
$1$ & $1$ & $1$ & $1$ & $1$ & $1$ \\ 
$2$ & $1$ & $5$ & $5$ & $3$ & $5$ \\ 
$3$ & $1$ & $5$ & $5$ & $2$ & $5$ \\ 
$4$ & $1$ & $2$ & $3$ & $4$ & $5$ \\ 
$5$ & $1$ & $5$ & $5$ & $5$ & $5$%
\end{tabular}
\end{center}

Let $\Gamma =\{1\}$ and define $S\times \Gamma \times S\rightarrow S$ by $%
x1y=xy$ for all $x,y\in S,$ then $S$ is a $\Gamma $-AG-groupoid.

Define an $IFS$ $A=(\mu _{A},\gamma _{A})$ of an AG-groupoid $S$ as follows: 
$\mu _{A}(1)=\mu _{A}(2)=\mu _{A}(3)=0.3,$ $\mu _{A}(4)=0.1$, $\mu
_{A}(5)=0.4,$ $\gamma _{A}(1)=0.2,$ $\gamma _{A}(2)=0.3,$ $\gamma
_{A}(3)=0.4,$ $\gamma _{A}(4)=0.5,$ $\gamma _{A}(5)=0.2.$ Now again define
an $IFS$ $B=(\mu _{B},\gamma _{B})$ of an AG-groupoid $S$ as follows: $\mu
_{B}(1)=\mu _{B}(2)=\mu _{B}(3)=0.5,$ $\mu _{B}(4)=0.4$, $\mu _{B}(5)=0.6,$ $%
\gamma _{B}(1)=0.3,$ $\gamma _{B}(2)=0.4,$ $\gamma _{B}(3)=0.5,$ $\gamma
_{B}(4)=0.6,$ $\gamma _{B}(5)=0.3.$ Then it is easy to observe that $A=(\mu
_{A},\gamma _{A})$ and $B=(\mu _{B},\gamma _{B})$ are an intuitionistic
fuzzy $\Gamma $-two-sided ideals of $S$ such that $(\mu _{A}\circ _{\Gamma
}\mu _{B})(a)=\{0.1,$ $0.3,$ $0.4\}=(\mu _{A}\cap \mu _{B})(a)$ for all $%
a\in S$ and similarly $(\gamma _{A}\circ _{\Gamma }\gamma _{B})(a)=(\gamma
_{A}\cap \gamma _{B})$ for all $a\in S$, that is, $A\circ _{\Gamma }B=A\cap
B $ but $S$ is not an intra-regular because $3\in S$ is not an intra-regular.

An $IFS$ $A=(\mu _{A},\gamma _{A})$ of an AG-groupoid is said to be
idempotent if $\mu _{A}\circ _{\Gamma }\mu _{A}=\mu _{A}$ and $\gamma
_{A}\circ _{\Gamma }\gamma _{A}=\gamma _{A},$ that is, $A\circ _{\Gamma }A=A$
or $A^{2}=A.$

\begin{lemma}
\label{idem}Every intuitionistic fuzzy $\Gamma $-two-sided ideal $A=(\mu
_{A},\gamma _{A})$ of an intra-regular $\Gamma $-AG-groupoid $S$ is $\Gamma $%
-idempotent.
\end{lemma}

\begin{proof}
Let $S$ be an intra-regular $\Gamma $-AG-groupoid and let $A=(\mu
_{A},\gamma _{A})$ be an intuitionistic fuzzy $\Gamma $-two-sided ideal of $%
S.$ Now for $a\in S$ there exists $x,y\in S$ and $\alpha ,\beta ,\delta \in
\Gamma $ such that $a=(x\alpha (a\beta a))\delta y.$ Now let $\xi \in \Gamma
,$ then by using $(3),$ $(5)$ and $(2),$ we have%
\begin{equation*}
a=(x\alpha (a\beta a))\delta y=(a\alpha (x\beta a))\delta (u\xi v)=(a\alpha
u)\delta ((x\beta a)\xi v).
\end{equation*}%
\begin{eqnarray*}
(\mu _{A}\circ _{\Gamma }\mu _{A})(a) &=&\dbigvee\limits_{a=(a\alpha
u)\delta ((x\beta a)\xi v)}\{\mu _{A}(a\alpha u)\wedge \mu _{A}((x\beta
a)\xi v)\}\geq \mu _{A}(a\alpha u)\wedge \mu _{A}((x\beta a)\xi v) \\
&\geq &\mu _{A}(a)\wedge \mu _{A}(a)=\mu _{A}(a).
\end{eqnarray*}

This shows that $\mu _{A}\circ _{\Gamma }\mu _{A}\supseteq \mu _{A}$ and by
using Lemma \ref{as}, $\mu _{A}\circ _{\Gamma }\mu _{A}\subseteq \mu _{A}$,
therefore $\mu _{A}\circ _{\Gamma }\mu _{A}=\mu _{A}.$ Similarly we can
prove that $\gamma _{A}\circ _{\Gamma }\gamma _{A}=\gamma _{A},$ which
implies that $A=(\mu _{A},\gamma _{A})$ is $\Gamma $-idempotent.
\end{proof}

\begin{theorem}
The set of intuitionistic fuzzy $\Gamma $-two-sided ideals of an
intra-regular $\Gamma $-AG$^{\ast \ast }$-groupoid $S$ forms a semilattice
structure with identity $\delta $, where $\delta =(S_{\delta },\Theta
_{\delta }).$
\end{theorem}

\begin{proof}
Let $\mathbb{I}_{\mu \gamma }$ be the set of intuitionistic fuzzy $\Gamma $%
-two-sided ideals of an intra-regular $\Gamma $-AG$^{\ast \ast }$-groupoid $%
S $ and let $A=(\mu _{A},\gamma _{A})$, $B=(\mu _{B},\gamma _{B})$ and $%
C=(\mu _{C},\gamma _{C})$ are any intuitionistic fuzzy $\Gamma $-two-sided
ideals of $\mathbb{I}_{\mu \gamma }.$ Clearly $\mathbb{I}_{\mu \gamma }$ is
closed and by Lemma \ref{idem}, we have $A\Gamma A=A$. Now by using Lemma %
\ref{fgh}, we get $A\circ _{\Gamma }B=B\circ _{\Gamma }A$ and therefore, we
have%
\begin{equation*}
(A\circ _{\Gamma }B)\circ _{\Gamma }C=(B\circ _{\Gamma }A)\circ _{\Gamma
}C=(C\circ _{\Gamma }A)\circ _{\Gamma }B=(A\circ _{\Gamma }C)\circ _{\Gamma
}B=(B\circ _{\Gamma }C)\circ _{\Gamma }A=A\circ _{\Gamma }(B\circ _{\Gamma
}C).
\end{equation*}

It is easy to see from Corollary \ref{cor} that $\delta $ is an identity in $%
\mathbb{I}_{\mu \gamma }.$
\end{proof}

\end{document}